\documentclass[11pt]{article}
\usepackage[english]{babel}
\usepackage[hmarginratio=3:2]{geometry}
\usepackage[titletoc]{appendix}
\usepackage{authblk}
\usepackage{graphicx}
\usepackage{amsmath}
\usepackage{amsthm}
\usepackage{amssymb}
\usepackage[mathscr]{eucal}
\usepackage{textcomp}
\usepackage[latin1]{inputenc}
\usepackage{latexsym}
\usepackage{hyperref}
\usepackage[dvipsnames]{xcolor}
\usepackage{framed}
\usepackage{tikz-cd}
\usepackage{ifsym}
\usepackage{enumerate}
\usepackage{sectsty}
\usepackage{chngcntr}

\numberwithin{equation}{section}

\numberwithin{subsection}{section}

\allowdisplaybreaks[1]

\newtheorem{theorem}{Theorem}[section]

\newtheorem{corollary}[theorem]{Corollary}
\newtheorem{proposition}[theorem]{Proposition}
\newtheorem{observation}[theorem]{Observation}
\newtheorem*{remark}{Remark}
\newtheorem*{definition}{Definition}

\newenvironment{enumeratea}    
{\begin{enumerate}[\upshape (a)]}
{\end{enumerate}}

\newenvironment{enumerate1}
{\begin{enumerate}[\upshape (1)]}
{\end{enumerate}}

 \newcommand\cB{\mathcal{B}}
\newcommand\cC{\mathcal{C}}

\newcommand\cI{\mathcal{I}}

\newcommand\cO{\mathcal{O}}

\newcommand\cU{\mathcal{U}} 
 \newcommand\cX{\mathcal{X}}
\newcommand\cY{\mathcal{Y}}

\newcommand\GG{\mathbb{G}}

\newcommand\QQ{\mathbb{Q}}

 \newcommand\ZZ{\mathbb{Z}}

\newcommand\rK{\mathrm{K}}

\newcommand\rS{\mathrm{S}}

\newcommand\rmm{\mathrm{m}}

 \newcommand\frm{\mathfrak{m}}

\newcommand\arr{\ifinner\to\else\longrightarrow\fi}

\newcommand\arrto{\ifinner\mapsto\else\longmapsto\fi}

\newcommand{\xarr}{\xrightarrow}

\newcommand{\eqdef}{\mathrel{\smash{\overset{\mathrm{\scriptscriptstyle def}} =}}}

\newcommand\into{\hookrightarrow}

\renewcommand{\th}{^\text{th}}

\def\displaytimes_#1{\mathrel{\mathop{\times}\limits_{#1}}}

\def\displayotimes_#1{\mathrel{\mathop{\bigotimes}\limits_{#1}}}

\newcommand\spec{\operatorname{Spec}}

\newcommand\rk{\operatorname{rk}}

\newcommand\pr{\operatorname{pr}}

\newcommand\double{\mathbin{\rightrightarrows}}

\newlength{\ignora}

\renewcommand{\setminus}{\smallsetminus}

\newcommand{\mmu}{\boldsymbol{\mu}}

\newcommand{\gm}{\GG_{\rmm}}

\newcommand{\GL}{\mathrm{GL}}

\newcommand{\abs}[1]{\left|#1\right|}

\newcommand{\step}[1]{\smallskip\textit{#1.}\;}

\renewcommand{\epsilon}{\varepsilon}

\renewcommand{\k}[1][*]{\operatorname{K}'_{#1}}

\newcommand{\kz}{\operatorname{K}_{0}}

\newcommand{\R}{\mathrm{R}}
\newcommand{\Rt}{\smash{\widetilde{\mathrm{R}}}}

\newcommand{\Q}[1]{\QQ(\zeta_{#1})}

\newcommand{\kmu}[2]{\rK'_{*}(#2)_{(\mmu_{#1})}}

\newcommand{\kg}[2][*]{\rK'_{#1}(#2)_{\bf g}}
\newcommand{\ka}[2][*]{\rK'_{#1}(#2)_{\bf a}}

\newcommand{\red}{_{\rm red}}

\title{The geometric K-theory of quotient stacks}
\author[1]{Angelo Vistoli}
\affil[1]{\textit{Scuola Normale Superiore, Piazza dei cavalieri 7, 56126, Pisa, Italy}}
\author[2]{Francesco Sala}
\affil[2]{\textit{Scuola Normale Superiore, Piazza dei cavalieri 7, 56126, Pisa, Italy}}
\author{Laurent Schadeck}
\date{}

\begin{document}
\maketitle

\begin{abstract}
Given a quotient of a regular noetherian separated algebraic space $X$ over a field by an affine algebraic group $G$ having finite stabilizers (with some mild technical conditions), G. Vezzosi and A. Vistoli defined the \emph{geometric part} of the rational equivariant K-theory $K(X,G)$ and conjectured that it is isomorphic to the rational K-theory of the quotient $X/G$. In this paper we refine the construction of geometric K-theory to the rational K-theory of a quotient stack $[X/G]$ over an arbitrary excellent base; we show that it is part of an intrinsic decomposition of the K-theory of the stack and prove many properties that make it amenable to computations.
\end{abstract}

\bigskip
\bigskip
\section*{Introduction}
In their paper \cite{vezzosi-vistoli02}, G. Vezzosi and A. Vistoli compared the equivariant K-theory of a separated noetherian algebraic space with that of its fixed loci by one-parameter diagonalizable subgroups.

They considered an action of an affine algebraic group $G$ on a regular separated noetherian algebraic space $X$ on a field $k$. Moreover they gave the following

\begin{definition}
The set of \emph{essential dual cyclic subgroups} $\cC(G)$ contains the conjugacy classes in $G$ of subgroups isomorphic to a group scheme of the form $\mmu_n:=Spec(k[x]/(x^n-1))$; a dual cyclic subgroup $\sigma$ is called \emph{essential} if $X^\sigma\neq\emptyset$.
\end{definition} 

Assuming a technical hypothesis (that the action is \emph{sufficiently rational}, see the introduction of \cite{vezzosi-vistoli02}), they proved a decomposition theorem for the rational equivariant K-theory $K(X,G)$: defining $\tilde R(\sigma)$ to be $\QQ(\zeta_n)$ when $\sigma$ has order $n$ and $w_G(\sigma):=N_G(\sigma)/C_G(\sigma)$ 
\[
K(X,G)\simeq \underset{\sigma\in \cC(G)}{\prod} (K(X^\sigma, C_G(\sigma))_{geom}\otimes \tilde R(\sigma))^{w_G(\sigma)}
\]
In the formula, $K(-)_{geom}$ is a suitable localization of the rational equivariant K-theory, which aims to substitute the rational K-theory of the geometric quotient space $(X/G)$.\\

In this paper we define the geometric rational K-theory of an algebraic stack $\cX$ with finite stabilizers, and we show that it is a piece in an intrinsic decomposition of the K-theory of the stack. Moreover given a presentation $\cX=[X/G]$, we compare this decomposition with the one induced by the dual cyclic subgroup schemes in $\cC(G)$.

\subsection*{Conventions and notations}

In what follows we will assume that $\cX$ is a separated algebraic stack with finite stabilizers, and that it can be written as a quotient $\cX=[X/G]$ where $X$ is a noetherian algebraic space over an affine connected base scheme $A=Spec(R)$, with $R$ excellent, and $G=GL_{n,R}$. 
The last condition is actually equivalent to the apparently weaker one that $\cX=[X/G]$, where $G$ is a flat $R-$linear algebraic group: indeed by choosing an embedding $G \subseteq \GL_{n} \eqdef \GL_{n,R}$, and replacing $X$ with $X\times^{G}\GL_{n}$, we may suppose that $G = \GL_{n}$. We will always write $\cX=[X/G]$, where $G$ is a product of general linear groups on $R$.

As in \cite{vezzosi-vistoli02}, we call a subgroup scheme $\sigma \subseteq G$ \emph{dual cyclic} when it is isomorphic to $\mmu_{r, R}$ for some $r \geq 1$, and we let $\cC(G)$ be the set of essential dual cyclic subgroup schemes (as above, essential means that $X^\sigma\neq\emptyset$).\\

We specify that a $\GL_{n,R}$-subgroup scheme $\sigma$ isomorphic to $Spec(R[x]/(X^n-1))$ is called \emph{constant} if the decomposition of $R^n$ in eigenspaces consists of free modules.

We will call $R(\sigma)$ the rational ring of characters of $\sigma$, and $R(G)$ the Grothendieck ring of highest weight representations of $G=\GL_{n,R}$ (which is isomorphic to $\ZZ[\hat{T}]^W$, where $\hat T$ denotes the characters of a maximal torus and $W$ is the Weyl group of $G$).\\

There is a special class of stacks and morphisms for which the properties of the decomposition discussed above are stronger: this is the class of \emph{tame} algebraic stacks: following \cite{dan-olsson-vistoli}, such a stack is called tame if the stabilizer of any point is linearly reductive. 

In general, a morphism of stacks $f\colon\cX\rightarrow \cY$ is called relatively tame if the relative inertia $\cI_{\cY}\cX$ has linearly reductive geometric fibers (see \cite{dan-olsson-vistoli}). For such a morphism, the push-forward $f_*$ descends to K-theory.

\section{The intrinsic decomposition of the K-theory\\of a quotient stack}

Fix now $G=\GL_{n,R}$.

For any $r \geq 1$ we have a decomposition $\R\mmu_{r} = \prod_{s|r}\Q{s}$. We denote by $\R \sigma\arr\Rt\sigma$ the projection onto the $\Q{r}$ factor.
Recall from \cite{vezzosi-vistoli02} that if $\sigma \subseteq G$ is a constant dual cyclic subgroup, $\frm_{\sigma} \subseteq \R G$ is the maximal ideal that is the kernel of the composite $\R G \arr \R\sigma \arr \Rt\sigma$; this only  depends on the \emph{type} of $\sigma$. 
If $M$ is an $\R G$-module, the \emph{$\sigma$-localization} $M_{\sigma}$ is the localization $M_{\frm_{\sigma}}$. If $\sigma \subseteq G$ is a constant dual cyclic subgroup, then $\k(X, G)_{\sigma} \neq 0$ if and only if $X^{\sigma} \neq \emptyset$.

We define the multiplicative system $\Sigma_{r}^{\cX} = \Sigma_{r} \subseteq \kz \cX$ as follows. An element $\alpha \in \kz\cX$ is in $\Sigma_{r}$ if for all representable morphisms $\phi\colon \cB_{S}\mmu_{r} \arr \cX$, where $S$ is an $R$-scheme, the projection of $\phi^{*}a \in \R\mmu_{r}$ in $\Rt\mmu_{r} = \Q{r}$ is non-zero.

In particular, $\Sigma_{1}$ admits the following description. If $\cX_{1}$, \dots,~$\cX_{m}$ are the connected components of $\cX$, then $\k\cX = \bigoplus_{i=1}^{m}\k\cX_{i}$. There is a rank map $\rk_{i}\colon \k\cX_{i} \arr \QQ$. Then $\alpha$ is in $\Sigma_{1}$ if and only if $\rk_{i}\alpha \neq 0$ for all $i$.

A morphism $f\colon \cY \arr \cX$ of stacks induces a pullback homomorphism $f^{*}\colon \kz\cX \arr \kz\cY$ (here by $\kz$ we denote the naive ring of locally free sheaves on $\cX$). If $f$ is representable, one immediately sees that $f^{*}$ carries $\Sigma_{r}^{\cX}$ into $\Sigma_{r}^{\cY}$. If $f$ is not representable this is certainly not true in general, but it is if $r = 1$. 

\begin{definition}
The $\mmu_{r}$-localization $\kmu r \cX$ of $\k\cX$ is the $\kz\cX$-module $\Sigma_{r}^{-1}\k(\cX)$.
\end{definition}

Set $\k(X,G)_{(r)} \eqdef \prod_{\sigma\in \cC_{r}(G)}\k(X, G)_{\sigma}$, so that $\k(X,G) = \prod_{r \geq 1} \k(X,G)_{(r)}$. The support of $\k(X, G)_{\sigma}$ in $\spec \R G$ is the maximal ideal $\frm_{\sigma} = \ker(\R G \arr \Rt\sigma)$.\\

Let $\sigma_{1}$, \dots,~$\sigma_{m}$ be (constant) representatives of dual cyclic subgroups $\sigma \in \cC(G)$  with $\abs{\sigma} = r$ such that $\k(X, G)_{\sigma} \neq 0$. Then
\[
   \k(X,G)_{(r)}=\prod_{\substack {\sigma\in \cC(G)\\ \abs{\sigma} = r}}\k(X, G)_{\sigma} =
   \prod_{i=1}^{m}\k(X, G)_{\sigma_{i}}.
\]

Suppose that $\phi\colon\cB_{S}\sigma \arr \cX$ is a representable morphism, where $\sigma$ is a constant dual cyclic group over $S$. We may assume that there is an embedding $\sigma \subseteq G_{S}$, and a rational point $p \in X(S)$ which is fixed under the action of $\sigma$. Hence $X_{S}^{\sigma} \neq \emptyset$, so $\sigma$ is conjugate to some $\sigma_{i}$ (that is it belongs to the same type). We say in this case that $\phi$ is \emph{associated} to $\sigma_i$.

Let $\sigma$ be any of the $\sigma_i$'s; we define another multiplicative system $\Sigma_\sigma^\cX=\Sigma_\sigma\subseteq \kz \cX$ as follows: $\alpha\in\kz\cX$ lies in $\Sigma_\sigma$ if and only if for any $S/A$ and any $\phi\colon \cB_S\mmu_r\arr \cX$ associated to $\sigma$ the pullback $\phi^*\alpha\in \R\mmu_r$ has a nonzero projection to $\Rt\mmu_r$.

\begin{proposition}\label{prop:decomposition}
Let $r$ be a fixed positive integer and $\sigma \in \cC(G)$ be a constant dual cyclic subgroup with $\abs{\sigma} = r$ such that $\k(X, G)_{\sigma} \neq 0$.
\begin{enumerate1}
\item The projection $\k(X, G) \arr \prod_{\sigma \in \cC(G)}\k(X, G)_{\sigma}$ is an isomorphism.
\item The $\mmu_{r}$-localization $\k \cX \arr \kmu r \cX$ isomorphic to the projection $\k(X,G) \arr \k(X,G)_{(r)}$.
\item The localization $\Sigma_\sigma^{-1}\k(\cX)$ is isomorphic to $\k(X,G)_\sigma$.
\end{enumerate1}

\end{proposition}

As a corollary, we get the following.

\begin{theorem}
The projections $\k(\cX) \arr \kmu r {\cX}$ induce an isomorphism
   \[
   \k(\cX) \simeq \prod_{r \geq 1}\kmu r {\cX}\,.
   \]
\end{theorem}

Another consequence is this.

\begin{corollary}\label{vanishing}
We have $\kmu r \cX \neq 0$ if and only if there exists a geometric point of $\cX$ whose automorphism group scheme contains a dual cyclic subgroup of order $r$.

In particular when $\cX$ is tame, $\kmu r \cX = 0$ for all $r > 1$ if and only if $\cX$ is an algebraic space.
\end{corollary}

\begin{proof}[Proof of Proposition~\ref{prop:decomposition}]
We are going to prove the three stataments jointly. Let $\sigma_{1}$, \dots,~$\sigma_{m}$ be all the constant dual cyclic subgroups $\sigma \in \cC(G)$  with $\abs{\sigma} = r$ such that $\k(X, G)_{\sigma} \neq 0$.

Call $S_\sigma:=\R G\setminus \frm_\sigma$ and $S_{r} \subseteq \R G$ the multiplicative system $\R G \setminus (\frm_{\sigma_{1}} \cup \dots\cup \frm_{\sigma_{m}})$; the induced map $S_{r}^{-1}\k(X, G) \arr \prod_{i=1}^{m}\k(X, G)_{\sigma_{i}}$ is an isomorphism (this follows easily from the fact that the support of $\k(X, G)$ in $\spec\R G$ consists of a finite number of closed points).

We claim that the images of $S_{r},S_\sigma \subseteq \R G$ in $\kz(X, G)$ through the homomorphism $\R G \arr \kz(X, G)$ are contained in $\Sigma_{r},\Sigma_\sigma$ respectively. In fact, let $\alpha \in S_{r}$, suppose that $\cB_{S}\sigma \arr \cX$ is a representable morphism, where $\sigma$ is a dual cyclic group over $S$; it is associated to one of the $\sigma_{i}$. The morphisms $\R G \arr \Rt \sigma_{i} = \Q{r}$ defined by the embedding $\sigma_{i} \subseteq G$ and the composite
   \[
   \R G \arr \kz(X,G) \arr \kz \cB_{S}\sigma = \R \sigma \arr \Rt\sigma = \Q{r}
   \]
defined by the morphism $\cB_{S}\sigma \arr \cX$ coincide; since by hypothesis $\alpha$ does not map to $0$ in $\Q{r}$, it follows that the image of $\alpha$ in $\kz\cX$ is in $\Sigma_{r}$, as claimed.

If $\alpha$ lies in $\Sigma_\sigma$ the reasoning is completely analogous.\\

Notice that if $\alpha \in \kz\cX = \kz(X, G)$, the multiplication by $\alpha$ on $\k(X, G)$ gives a homomorphism of $\R G$-modules, so that it preserves the multiplication above. We get a factorization
   \[
   \begin{tikzcd}[row sep = 10pt]
   \k(X,G) \ar[r] & S_{r}^{-1}\k(X,G) \ar[r]\ar[d, equal]&
   \Sigma_{r}^{-1}\k(X, G)\ar[d, equal]\\
   & \displaystyle\prod_{\substack {\sigma\in \cC(G)\\ \abs{\sigma} = r}}\k(X, G)_{\sigma} &
   \kmu r {X,G}
   \end{tikzcd}
   \]
we need to show that the resulting homomorphism
   \[
   \prod_{\substack {\sigma\in \cC(G)\\ \abs{\sigma} = r}}\k(X, G)_{\sigma}
   \arr \kmu r {X,G}
   \]
is an isomorphism.

Similarly we have a morphism $ \k(X,G)_\sigma\arr\Sigma_\sigma^{-1}\kz(\cX)$ and we want it to be a bijection.\\

These are equivalent to the following. First of all, notice that if $\alpha \in \kz(X, G)$, multiplication by $\alpha$ gives an endomorphism of the $\R G$-modules $\k(X, G)$, hence it descends to endomorphisms of $\k(X, G)_{\sigma}$ for each $\sigma \in \cC(G)$. We need to prove that the if $\alpha$ is in $\Sigma_{r}$ or in $\Sigma_\sigma$, then $\alpha$ induces an automorphism of $\k(X, G)_{\sigma}$. Actually, $\Sigma_r\subseteq\Sigma_\sigma$, so me may just suppose $\alpha\in\Sigma_\sigma$.

The proof of this fact is fairly long, and will be split into steps.

\step{Step 1} Assume that $G$ is a finite diagonalizable group over $R$ acting trivially on an affine scheme $X$. Then $\cX$ is of the form $X \times_A \cB_{R}G$ (here the tensor products are intended over $K_0(A)$). Then we have $\k\cX = \k X \otimes \R G$ and $\kz(\cX) = \kz X \otimes \R G$. Then $\cC(G)$ is the set of dual cyclic subgroups of $G$. We have a decomposition $\R G = \prod_{\sigma \in \cC(G)} \Rt\sigma$. Then $\k(X, G)_{\sigma} = \k X\otimes\Rt\sigma$ and $\k(X,G)_{(r)} = \k X\otimes \prod_{\substack {\sigma\in \cC(G)\\ \abs{\sigma} = r}}\Rt\sigma$. In particular, $(1)$ follows immediately.

The action of $\kz(X,G)$ on $\k(X, G)_{\sigma}$ is induced by the action of $\kz(X,G) = \kz X \otimes \R\sigma$ on $\k(X, G) = \k X \otimes \Rt\sigma$, and this in turn is induced by the action of $\kz X$ on $\k X$. This factors through an action of $\kz X \otimes \Rt\sigma$; thus it is enough to show that the image of $\alpha$ in $\kz X \otimes \Rt\sigma$ is invertible.

Let $X_{1}$, \dots,~$X_{m}$ be the connected components of $X$. Then $\kz X = \prod_{i=1}^{m}\kz X_{i}$. There is a rank homomorphism $\rk\colon \kz X_{i} \arr \QQ$, whose kernel is nilpotent (\cite{Fulton-Lang}, V, sect. 3); from this we obtain a homomorphism
   \[
   \k X\otimes \Rt\sigma \arr (\Rt\sigma)^{m}
   \]
with nilpotent kernel. Hence it is enough to show that the image of $\alpha$ in each copy of $\Rt\sigma$ is non-zero. But $\alpha_{i}$ is obtained as follows: choose a point $\spec S \arr X_{i}$. This gives a morphism $\cB_{S}G_{S} \arr X \times \cB_{R}G$; by composing with the morphism $\cB_{S}\sigma_{S} \arr \cB_{S}G_{S}$ induced by the embedding $\sigma \subseteq G$; the resulting homomorphism $\kz(X,G) \arr \Rt\sigma$ is immediately seen to coincide with the $i\th$ homomorphism $\k X\otimes \Rt\sigma \arr \Rt\sigma$. Since this is not zero, by hypothesis, this concludes the proof.

\step{Step 2: the case of a torus action} Here we assume that $G = \gm^{n}$ is a split torus. If $Y \subseteq X$ is a $G$-invariant subscheme of $X$, we denote by $\alpha_{Y}$ the restriction of $\alpha$ to $\kz(Y,G)$. By noetherian induction, we can assume that $\alpha_{Y}$ induces an automorphism of $\k(Y, G)_{\sigma}$ for all proper closed $G$-invariant subschemes $Y \arr X$. By \cite{Thomason} there exists an open $G$-invariant subscheme $U$ that is $G$-equivariantly isomorphic to a scheme of the form $V \times (G/\Gamma)$, where $\Gamma \subseteq G$ is a finite diagonalizable subgroup scheme, $V$ is a scheme over $R$, and the action of $G$ on $V \times (G/\Gamma)$ is induced by the trivial action of $G$ on $V$ and the action on $G/\Gamma$ by translation. By shrinking again $V$ we may assume that $V$ is affine and connected. 

Assume that $X = U$. In this case $\k(X, G) = \k(U, \Gamma) = \k U \otimes \R\Gamma$. If $\sigma \nsubseteq \Gamma$, then $X^{\sigma} =  \emptyset$, so that $\k(X,G)_{\sigma} = 0$, and the result is obvious. If $\sigma \subseteq \Gamma$, then it is easy to convince oneself that $\k(X,G)_{\sigma} = \k(V, \Gamma)_{\sigma}$, and so the result follows from the preceding case. 

If $X \neq U$, call $Y$ the complement of $U$, with its reduced scheme structure. This is $G$-invariant. For each $i \geq 0$ we get a commutative diagram
   \[
   \begin{tikzcd}[column sep= 17pt]
   \k[i+1](U, G) \ar[r]\ar[d]&
      \k[i](Y,G) \ar[r]\ar[d]
      & \k[i](X, G) \ar[r]\ar[d]
      & \k[i](U, G) \ar[r]\ar[d]
      &\k[i-1](Y, G)\ar[d]\\
\prod_\sigma  \k[i+1](U, G)_\sigma \ar[r]& \prod_\sigma\k[i](Y,G)_\sigma \ar[r]& \prod_\sigma \k[i](X, G)_\sigma \ar[r]
      & \prod_\sigma\k[i](U, G)_\sigma \ar[r] & \prod_\sigma\k[i-1](Y, G)_\sigma
   \end{tikzcd}
   \]
with exact rows. Since the result is true for $Y$ and $U$, so that the first, second, fourth and fifth columns are isomorphisms, we see that the central arrow is also an isomorphism by the "five lemma"; so we have $(1)$ as claimed.\\

As for the remaining points, we also have a diagram
   \[
   \begin{tikzcd}[column sep= 17pt]
   \k[i+1](U, G) \ar[r]\ar[d, "\cdot\alpha_{U}"]&
      \k[i](Y,G) \ar[r]\ar[d, "\cdot\alpha_{Y}"]
      & \k[i](X, G) \ar[r]\ar[d, "\cdot\alpha"]
      & \k[i](U, G) \ar[r]\ar[d, "\cdot\alpha_{U}"]
      &\k[i-1](Y, G)\ar[d, "\cdot\alpha_{Y}"]\\
   \k[i+1](U, G) \ar[r]& \k[i](Y,G) \ar[r]& \k[i](X, G) \ar[r]
      & \k[i](U, G) \ar[r] &\k[i-1](Y, G)
   \end{tikzcd}
   \]
with exact rows. As before, the results follow.

\step{Step 3: the case $G = \GL_{n}$} Consider the standard maximal torus $T \subseteq G$ with its Weyl group $\rS_{n}$. The subgroup $\sigma \subseteq G$ is conjugate to a subgroup of $T$, well defined, up to conjugation by an element of $\rS_{n}$ thanks to the following easy

\begin{observation}\label{lemma:centralizer}
$\cC(G)$ is in natural bijective correspondence with the set of orbits for the action of $\rS_{n}$ on $\cC(T)$.
\end{observation}

 We have that the restriction map induces an isomorphism $\k(X, G) = \k(X,T)^{\rS_{n}}$ (see Appendix B); moreover the map $\R\GL_{n}\arr \R T$ is a Galois cover with Galois group $S_n$. In particular for any $\R T$-module with finite support and $\sigma\in\cC(\GL_{n})$ there is an isomorphism $M_\sigma\simeq\underset{\sigma '\in S_n\sigma}{\prod}M_{\sigma '}$.

We conclude that
   \[
   \prod_{\sigma \in \cC(G)} \k(X,G)_{\sigma} \simeq
   \bigl(\prod_{\sigma \in \cC(T)}\k(X, T)_{\sigma}\bigr)^{\rS_{n}}\,.
   \]
and $(1)$ follows.

If $\Gamma \subseteq \rS_{n}$ is the stabilizer of $\sigma \subseteq T$ under the action of $\rS_{n}$, we deduce that
   \[
   \k(X, G)_{\sigma} \simeq \k(X, T)_{\sigma}^{\Gamma}\,.
   \]
If $\beta$ is the image in $\kz(X, T)$ of $\alpha \in \kz(X, G)$, then multiplication by $\beta$ on $\k(X,T)_{\sigma}$ is $\Gamma$-equivariant, and the its restriction to $\k(X, G)_{\sigma} = \k(X, T)_{\sigma}^{\Gamma}$ is multiplication by $\alpha$. Since multiplication by $\beta$ on $\k(X,T)_{\sigma}$ is an automorphism, by the previous step, we deduced that multiplication by $\alpha$ is also an automorphism, as claimed.

\end{proof}

In the decomposition $\k\cX = \prod_{r}\kmu r \cX$ we split off the factor $\kmu 1 \cX$, which we call the \emph{geometric K-theory} of $\cX$, and denote by $\kg\cX$, and the factor $\prod_{r \geq 2}\kmu r \cX$, which we call the \emph{algebraic part} of the K-theory of $\cX$, and denote by $\ka\cX$. Thus we have the \emph{fundamental decomposition}
   \[
   \k\cX = \kg\cX \oplus \ka\cX\,.
   \]
As we already remarked, $\ka\cX = 0$ if and only if $\cX$ is an algebraic space.

\begin{proposition}\label{prop:functoriality-fundamental-decomposition}
Let $f\colon \cY \arr \cX$ be a representable morphism of stacks.

\begin{enumerate1}

\item Assume that $f$ has finite flat dimension. Then the pullback $f^{*}\colon \k\cX \arr \k\cY$ preserves the fundamental decomposition.

\item Assume that $f$ is proper. Then the pushforward $f_{*}\colon \k\cY \arr \k\cX$ preserves the fundamental decomposition.

\end{enumerate1}
\end{proposition}

\begin{proof}
Write $\cX = [X/G]$. Then we have the decomposition
   \[
   \k\cX = \k(X, G) = \prod_{\sigma \in \cC(G)}\k(X, G)_{\sigma}\,;
   \]
we have $\kg\cX = \k(X, G)_{\{1\}}$ and $\ka\cX = \prod_{\substack{\sigma \in \cC(G)\\ \sigma \neq \{1\}}}\k(X, G)_{\sigma}$; hence it is enough to show that $f^{*}$ and $f_{*}$ preserve the decomposition
   \[
   \k(X, G) = \prod_{\sigma \in \cC(G)}\k(X, G)_{\sigma}\,.
   \]

Set $Y \eqdef \cY\times_{\cX}X$; then $Y$ is an algebraic space, since $f$ is representable; furthermore there is an action of $G$ on $Y$ such that $\cY = [Y/G]$ and the map $f$ is $G$-equivariant. Then $\k\cY = \k(Y,G)$ and $\k\cY$ acquires a structure of $\R G$-module. Furthermore, $f^{*}$ and $f_{*}$ are homomorphism of $\R G$-modules, and the result follows from this.
\end{proof}

Thus, if $f\colon \cY \arr \cX$ is representable of finite flat dimension, this defines homomorphisms  $f^{*}\colon \kg\cX \arr \kg\cY$ and $f^{*}\colon \ka\cX \arr \ka\cY$, while if $f$ is proper it defines homomorphism $f_{*}\colon \kg\cY \arr \kg\cX$ and $f_{*}\colon \ka\cY \arr \ka\cX$. These make $\kg-$ and $\ka-$ into contravariant functors for representable maps of finite flat dimension, and covariant functors for proper representable maps.\\

If $f$ is not representable, it is immediate to give examples to show that the result above can fail. However, we have the following important fact. The fundamental decomposition $\k\cX = \kg\cX \oplus \ka\cX$ gives both a projection $\k\cX\arr\kg\cX$ and an embedding $\kg\cX \arr \k\cX$. 

\begin{proposition}\label{prop:functoriality}
Let $f\colon \cY \arr \cX$ be a homomorphism of stacks.

\begin{enumerate1}

\item Assume that $f$ has finite flat dimension. Then there exists a homomorphism
   \[
   f^{*}\colon \kg\cX \arr \kg\cY
   \]
such that the diagram
   \[
   \begin{tikzcd}
   \k\cX \ar[r, two heads]\ar[d, "f^{*}"]& \kg\cX\ar[d, "f^{*}"]\\
   \k\cY \ar[r, two heads] & \kg\cY
   \end{tikzcd}
   \]
commutes.

\item Suppose that $f$ is proper and relatively tame. Then there exists a homomorphism
   \[
   f_{*}\colon \kg\cY \arr \kg\cX
   \]
such that the diagram
   \[
   \begin{tikzcd}
   \kg\cY \ar[r, hook]\ar[d, "f_{*}"]& \k\cY\ar[d, "f_{*}"]\\
   \kg\cX \ar[r, hook] & \k\cX
   \end{tikzcd}
   \]
commutes.

The maps $\kg- \into \k(-)$ and\/ $\k(-) \twoheadrightarrow \kg-$  in the diagrams above are the embeddings and the  projections coming from the fundamental decomposition. 
\end{enumerate1}
\end{proposition}

The homomorphism $f^{*}$ and $f_{*}$ defined above are clearly unique; they make $\kg-$ into a contravariant functor for maps of finite flat dimension, and a covariant functor for proper maps.

\begin{proof}
The first part is a consequence of the fact that $f\colon \kz\cX \arr \kz\cY$ carries $\Sigma_{1}^{\cX}$ into $\Sigma_{1}^{\cY}$.

For the second part, write $\cX = [X/G]$ and $\cY = [Y/H]$, where $X$ and $Y$ are schemes, and $G$ and $H$ are affine algebraic groups acting sufficiently rationally on $X$ and $Y$. The projection $X \arr \cX$ and the composite $Y \arr \cY \xarr{f} \cX$ are respectively $G$ and $H$ invariant. Set $Z = X \times_{\cX}Y$; there is a natural action of $G \times H$ on $Z$, and it is easy to see that $[Z/G\times H] = \cY$. The projection $Z \arr X$ is equivariant for the projection $\pr_{1}\colon G \times H \arr G$, and the induced morphism
   \[
   \cY = [Z/ G\times H] \arr [X/G] = \cX
   \]
is isomorphic to $f$. The homomorphism $f_{*}\colon \k(Z, G \times H) \arr \k(X, G)$ is a homomorphism of $\R G$-modules, where $\k(Z, G \times H)$ is considered as an $\R G$-module through the homomorphism $\pr_{1}^{*}\colon \R G \arr \R(G \times H)$.

Consider the decompositions
   \[
   \k(Z, G \times H) = \prod_{\rho \in \cC(G \times H)}\k(Z, G \times H)_{\rho}
   \]
and
   \[
   \k(X, G) = \prod_{\sigma \in \cC(G \times H)}\k(X, G)_{\sigma}\,.
   \]
If $\sigma \in \cC(G)$, we can also consider the $\sigma$-localization
   \[
   \k(Z, G \times H)_{\sigma} = (\R G \setminus \frm_{\sigma})^{-1}\k(Z, G \times H)
   \]
of $\k(Z, G \times H)$; it is immediate to see that it coincides with the quotient
   \[
   \prod_{\substack{\rho \in \cC(G\times H)\\
   \rho \twoheadrightarrow \sigma}}\k(Z, G \times H)_{\rho}
   \]
of $\k(Z, G \times H)$. If
   \[
   \eta \in \kg\cY = \k(Z,G \times H)_{\{1\}} \subseteq \k(Z, G \times H)\,
   \]
the by definition the image of $\eta$ in $\k(Z, G \times H)_{\sigma}$ is zero for every $\sigma \in \cC(G)$ with $\sigma \neq \{1\}$; this implies that the image of  $f_{*}\eta \in \k\cX$ in 
   \[
   \ka\cX = \prod_{\substack{\rho \in \cC(G \times H)\\ \rho \neq \{1\}}}\k(X,G)_\sigma
   \]
 is zero, which implies $f_{*}\xi \in \kg\cX$, as claimed.
\end{proof}

\begin{remark}
\label{remark:push-forward}
Actually, we can make the above proposition more precise. Suppose that we have a morphism of stacks $f\colon [Z/G\times H]\arr [X/G]$ induced by an equivariant map $Z\arr X$.

Let $\rho$ be a dual cyclic subgroup of $G\times H$ and $\sigma$ be its projection to $G$. 

\begin{enumerate1}
\item The pull-back $f^*$ is a map of $\R G$-modules, so it induces a map
\[
f^*:\k(X,G)_\sigma=S_\sigma^{-1}\k(X,G)\arr S_\sigma^{-1}\k(Z,G\times H)=\prod_{\substack{\rho \in \cC(G\times H)\\
   \rho \twoheadrightarrow \sigma}}\k(Z, G \times H)_{\rho}
\]
\item The push-forward $f_*\colon \k(Z,G\times H)\arr \k(X,G)$ is a morphism of $\R G$-modules. If $\eta\in \k(Z, G\times H)$ has support $\rho$ (so that its projections to the $\rho '$-localizations are zero for any $\rho '\neq\rho$), then its image $f_*(\eta)$ projects to zero in $\k(X,G)_{\sigma '}$ for every $\sigma '\neq\sigma$: indeed for any such $\sigma '$ and $\rho '\twoheadrightarrow \sigma '$ it holds that $\rho '\neq \rho$, and $\eta$ has null image in $\k(Z,  G\times H)_{\rho '}$. We conclude that $\eta$ is zero when localized at $\sigma '$, since $\k(Z,G\times H)_{\sigma '}=\prod_{\substack{\rho ' \in \cC(G\times H)\\
   \rho ' \twoheadrightarrow \sigma '}}\k(Z, G \times H)_{\rho '}$; being $f_*$ a map of $\R G-$modules, $f_*(\eta)$ is also zero when localized at $\sigma '$.

In particular $f_*\eta$ has support at $\sigma$, and thus lies in $\k(X,G)_\sigma$.
\end{enumerate1}

\end{remark}

We will now prove some following facts.

\begin{proposition}\label{prop:reduced-pushforward}
Let $i\colon \cX\red \into\cX$ be the reduction of $\cX$. Then
   \[
   i_{*}\colon \kg{\cX\red} \arr \kg\cX
   \]
is an isomorphism.
\end{proposition}

\begin{proof}
This follows from the fact that, by dévissage, $i_{*}\colon \k{\cX\red} \arr \k\cX$ is an isomorphism.
\end{proof}

\begin{proposition}\label{prop:descent-geometric-K}
Let $\pi\colon \cX' \arr \cX$ a representable finite faithfully flat morphism of  stacks. Then the sequences
   \[
   0 \arr \kg\cX \xarr{\pi^{*}} \kg{\cX'}
   \xarr{\pr^{*}_{1} - \pr^{*}_{2}}
   \kg{\cX'\times_{\cX} \cX'}   
   \]
and
   \[
   \kg{\cX'\times_{\cX} \cX'}
   \xarr{\pr_{1*} - \pr_{*2}} \kg{\cX'}
   \xarr{\pi_{*}}
   \kg\cX \arr 0   
   \]
are exact.
\end{proposition}

This is easily seen to fail for K-theory, for example, for the representable finite faithfully flat morphism $\spec k \arr \cB_{k}G$, where $G$ is a nontrivial finite group scheme.

\begin{corollary}\label{finite-cover}
Let $\Gamma$ be a finite group, $\pi\colon \cX'\arr \cX$ a Galois cover with group $\Gamma$. Then the pullback
   \[
   \pi^{*} \colon \kg\cX \arr \kg{\cX'}^{\Gamma}
   \]
and the pushforward
   \[
   \pi_{*}\colon (\kg{\cX'})_{\Gamma} \arr \kg\cX
   \]
are isomorphisms.
\end{corollary}

Here of course $\kg{\cX'}^{\Gamma}$ is the group of invariants, and $\pi_{*}(\kg{\cX'})_{\Gamma}$ is the group of covariants.

This applies in particular when $\cX$ is a quotient stack of the form $[X/\Gamma]$, where $\Gamma$ is a finite group acting on a scheme $X$. In this case we have isomorphism $\kg{[X/\Gamma]} \simeq (\k X)^{\Gamma}$ and $(\k X)_{\Gamma} \simeq \kg\cX$.

\begin{proof}[Proof of Proposition~\ref{prop:descent-geometric-K}]
Consider the class $\alpha \eqdef [\pi_{*}\cO_{\cX'}] \in \kz{\cX}$; because of the hypotheses, we have $\alpha \in \Sigma_{1}^{\cX}$. Hence, multiplication by $\alpha$ gives an automorphism of $\kg\cX$.

Suppose that $\xi\in \k\cX$, and $\pi^{*}\xi = 0$. Then by the projection formula we have $0 = \pi_{*}\pi^{*}\xi = \alpha\xi$; hence $\xi = 0$. This proves that $\pi^{*}$ is injective.

Now take $\xi' \in \k\cX'$ such that $\pr_{1}^{*}\xi' = \pr_{2}^{*}\xi'$. Then we have
   \begin{align*}
   \pi^{*}\pi_{*}\xi' &= \pr_{1*}\pr_{2}^{*}\xi'\\
   &= \pr_{1*}\pr_{1}^{*}\xi'\\
   &= [\pr_{1*}\cO_{\cX'\times_{\cX}\cO_{\cX}}]\xi'\\
   &= (\pi^{*}\alpha)\xi'.
   \end{align*}
Call $\eta \in \kg\cX$ the element such that $\alpha\eta = \pi_{*}\xi'$. Then we have
   \begin{align*}
   (\pi^{*}\alpha)\pi^{*}\eta &= \pi^{*}(\alpha\eta)\\
   &= \pi^{*}\pi_{*}\xi'\\
   &= (\pi^{*}\alpha)\xi'\,.
   \end{align*}
But we have $\pi^{*}\alpha \in \Sigma_{1}^{\cX'}$, hence $\xi' = \pi^{*}\eta$. This ends the proof of the first part.

For the second part, consider a class $\xi \in \kg\cX$. Then $\alpha\xi = \pi_{*}\pi^{*}\xi$. If $\eta \in \kg{\cX'}$ is such that $(\pi^{*}\alpha)\eta = \pi^{*}\xi$, then $\alpha\xi = \pi_{*}\bigl((\pi^{*}\alpha)\eta\bigr) = \alpha\pi_{*}\eta$, so $\xi = \pi_{*}\eta$, and $\pi_{*}$ is surjective.

Now take $\xi' \in \kg{\cX'}$ such that $\pi_{*}\xi' = 0$. Then $\pr_{2*}\pr_{1}^{*}\xi' = \pi^{*}\pi_{*}\xi' = 0$, while $\pr_{1*}\pr_{1}^{*}\xi' = (\pi^{*}\alpha)\xi'$. Denote by $\rho\colon \cX'\times_{\cX} \cX' \arr \cX$ is the composite morphism. If $\eta\in \kg{\cX'\tau_{\cX}\cX'}$ is such that $\rho^{*}\eta = \pr_{1}^{*}\xi'$ then we have $\pr_{1*}\eta = \xi'$, while $(\pi^{*}\alpha)\pr_{2*}\pr_{1}^{*}\eta = 0$, so $\pr_{2*}\pr_{1}^{*}\eta = 0$. Hence $(\pr_{1*} - \pr_{2*})\eta = \xi'$; this finishes the proof.
\end{proof}

Another way of thinking about $\kg{\cX} \subseteq \k\cX$ and $\ka{\cX} \subseteq \k\cX$ is the following: $\kg{\cX}$ is formed of elements that come from schemes via pushforward, while $\ka{X}$ is formed of all the elements that pull back to $0$ to schemes. This is probably true in general, but we can only prove it for tame with quasi-projective moduli space.

\begin{proposition}
Assume that a quotient stack $\cX$ is of finite type over $S$. Let $\xi$ be an element of $\k\cX$. Then
\begin{enumeratea}

\item $\xi \in \ka\cX$ if and only if for every morphism $f\colon Y \arr \cX$ of finite flat dimension, where $Y$ is an algebraic space, we have $f^{*}\xi = 0$, and

\item $\xi \in \kg\cX$ if and only if there exists a proper morphism $f\colon Y \arr \cX$, where $Y$ is an algebraic space, and a class $\eta\in \k Y$ such that $\xi = f_{*}\eta$.

\end{enumeratea}
\end{proposition}

\begin{proof}
Both ``only if'' directions are clear from the fact that $\ka Y = 0$, and from Proposition~\ref{prop:functoriality}. For the other direction we need the following fact.

\begin{proposition}[\hbox{\cite[Lemma~2.12]{KV}}]
Assume that $\cX$ is of finite type over $S$. Then there exists a flat finite surjective morphism $\pi\colon Y \arr \cX$, where $Y$ is an algebraic space over $S$.
\end{proposition}

If $\pi\colon Y \arr \cX$ is as above, then the pullback $\pi^{*}\colon \kg{\cX} \arr \k Y$ is injective, while the pushforward $\pi_{*}\colon \k Y \arr \kg\cX$ is surjective, by Proposition~\ref{prop:descent-geometric-K}. This proves what we want.
\end{proof}

Let $\pi\colon \cX \arr M$ be the moduli space of $\cX$; then we get a homomorphism
   \[
   \pi_{*}\colon \kg\cX \arr \kg M = \k M\,.
   \]

\begin{theorem}\label{thm:pushforward-moduli}
The pushforward $\pi_{*}\colon \kg\cX \arr \k M$ is an isomorphism when $\cX$ is tame.
\end{theorem}

\begin{proof}
Let $\cY \subseteq \cX$ be a closed substack. If $\cY \arr N$ is the moduli space of $\cY$, the natural morphism $N \arr M$ is a closed embedding, because $\cX$ is tame. By noetherian induction, we may assume that for every proper closed substack $\cY$, the proper pushforward $\kg\cY \arr \kg N$ is an isomorphism.

Consider the associated reduced substack $\cX\red \subseteq \cX$; its moduli space is $M\red$. We have a commutative diagram of proper pushforwards
   \[
   \begin{tikzcd}
   \kg{\cX\red} \ar[r]\ar[d]& \kg\cX \ar[d]\\
   \k{M\red} \ar[r]& \k M
   \end{tikzcd}
   \]
in which the rows are isomorphisms. If $\cX$ is not reduced then $\cX\red \neq\cX$, and the left hand column  is also an isomorphism; it follows that the result holds. So we may assume that $\cX$ is reduced.

Let $N \subseteq M$ be a closed subscheme, and set $U = M \setminus N$. Call $\cY$ and $\cU$ respectively the inverse images of $N$ and $U$ in $\cX$; then $\cY \arr N$ and $\cU \arr U$ are moduli spaces, by \cite{dan-olsson-vistoli}, 3.3. By noetherian induction we may assume that for every $N \subsetneqq M$, the pushforward $\kg{\cY} \arr \k N$ is an isomorphism. We have a commutative diagram with exact rows
   \[
   \begin{tikzcd}
   \kg[i-1]{\cU} \dar\rar
      & \kg[i]\cY \ar[r]\ar[d] & \kg[i]\cX \ar[r]\ar[d] & \kg[i]\cU \ar[r]\ar[d] &
    \kg[i-1]\cY \ar[d]\\
   \k[i-1]U \rar &\k[i]N \ar[r] & \k[i]M \ar[r] & \k[i]U \ar[r] &
    \k[i-1]N\\   
   \end{tikzcd}
   \]
in which the rows are parts of homotopy exact sequences, and the columns are given by proper pushforwards. If we assume that $U \neq \emptyset$ and that $\k\cU \arr \k U$ is an isomorphism, then the thesis follows. Thus in order the prove the thesis we may restrict $M$ to an non-empty open subset.

Now assume that $M' \arr M$ is an fppf cover, and set $M'' \eqdef M'\times_{M}M'$, $\cX' \eqdef M'\times_{M}\cX$, $\cX'' \eqdef M''\times_{M}\cX = \cX'\times_{\cX}\cX'$. We have a commutative diagram
   \[
   \begin{tikzcd}
   \kg\cX \ar[r]\ar[d] & \kg{\cX'} \ar[d]
   \arrow[r, shift left, "\pr_{1}^{*}"]
   \arrow[r, shift right, "\pr_{2}^{*}" below]
   & \kg{\cX''}\ar[d]\\
   \k M \ar[r] & \k{M'}
   \arrow[r, shift left, "\pr_{1}^{*}"]
   \arrow[r, shift right, "\pr_{2}^{*}" below]& \k{M''}\\
   \end{tikzcd}
   \]
in which the columns are push-forwards and the rows are equalizers, by Proposition~\ref{prop:descent-geometric-K}. It follows that if the result holds for $\cX'$ and $\cX''$, it also follows for $\cX$.

By \cite{KV}, Theorem 2.7, there exists a scheme $X$ and a finite morphism $X \arr \cX$; since $\cX$ is reduced, after passing to a nonempty open substack we may assume that $X \arr \cX$ and $\cX \arr M$ are flat. After passing to the finite fppf cover $X \arr M$, we can also assume that $\cX \arr M$ has a flat section $\sigma\colon M \arr \cX$. If $\Delta = M \times_{\cX}M \arr M$ denotes the automorphism group scheme of $\sigma$, then $\Delta \arr M$ is a linearly reductive finite group scheme over $M$. Since $\cX$ is the fppf quotient stack of the groupoid $M \times_{\cX}M \double M$, it is equivalent to the classifying stack $\cB_{M}\Delta \arr M$. By taking a restriction and a finite fppf cover, we may assume that $\Delta$ is of the form $\Delta_{1} \ltimes \Delta_{0}$, where $\Delta_{0} \arr M$ is finite and diagonalizable, and $\Delta_{1} \arr M$ is a tame étale constant group scheme (see \cite{dan-olsson-vistoli}); in the terminology of \cite{dan-olsson-vistoli}, such group schemes are called \emph{well-split}. Such a group scheme is the pullback of a well split group scheme $\Gamma \arr \spec k$.

So we are reduced to the case $\cX = M \times \cB_{S}\Gamma$, where $\Gamma \arr S$ is well-split. Then we have $\kg\cX = \k M \otimes (\R \Gamma)_{\{1\}}$, and the result follows if we show that the homomorphism $(\R \Gamma)_{\{1\}} \arr \QQ$ induced by the rank map $\R\Gamma \arr \QQ$ is an isomorphism. But from Proposition~\ref{prop:descent-geometric-K} applied to the morphism $S \arr \cB_{S}\Gamma$ we see that the pullback $(\R \Gamma)_{\{1\}} = \kg{\cB_{S}\Gamma} \arr \kg{S} $ is injective, which completes the proof.
\end{proof}

\end{document}